\documentclass[10pt]{IEEEtran4PSCC/IEEEtran4PSCC}

\ifCLASSINFOpdf
\else
\fi

\usepackage[cmex10]{amsmath}
\usepackage{psfrag}
\usepackage{multicol}
\usepackage{epstopdf}
\usepackage{graphicx,bm,upgreek,eucal}
\PassOptionsToPackage{bookmarks=false}{hyperref}
\usepackage{subfigure,textcomp,gensymb,url,multirow}
\usepackage{color,cite}
\usepackage{algorithm}
\usepackage{algpseudocode}
\usepackage{anysize, booktabs, setspace}
\marginsize{0.7in}{0.7in}{0.7in}{0.7in}
\usepackage{newpxmath}
\usepackage{amsfonts, amssymb}
\usepackage{cleveref, nomencl, blindtext, mathtools, etoolbox}
\usepackage[export]{adjustbox}

\newcommand{\comments}[1]{}
\newcommand{\problem}{MDOP}
\allowdisplaybreaks
\usepackage{tikz}
\usetikzlibrary{shapes,arrows,spy,positioning,patterns,shapes.geometric}
\usetikzlibrary{external}
\makeatletter
\def\munderbar#1{\underline{\sbox\tw@{$#1$}\dp\tw@\z@\box\tw@}}
\makeatother

\begin{document}
\title{Hierarchical Predictive Control Algorithms for \\ Optimal Design and Operation of Microgrids}
\author{
   Sai Krishna Kanth Hari$^{\dag}$, Kaarthik Sundar$^*$, Harsha Nagarajan$^\ddag$, Russell Bent$^*$, Scott Backhaus$^*$ \\
$^\dag$Department of Mechanical Engineering, Texas A\&M University, College Station, USA.\\
$^\ddag$Theoretical Division (T-5), Los Alamos National Laboratory, NM, USA.  Contact: harsha@lanl.gov \\
$^*$Center for Nonlinear Studies, Los Alamos National Laboratory, NM, USA.}

\maketitle

\begin{abstract}
In recent years, microgrids, \emph{i.e.}, disconnected distribution systems, have received increasing interest from power system utilities to support the economic and resiliency posture of their systems.  The economics of long distance transmission lines prevent many remote communities from connecting to bulk transmission systems and these communities rely on off-grid microgrid technology.  Furthermore, communities that are connected to the bulk transmission system are investigating microgrid technologies that will support their ability to disconnect and operate independently during extreme events. In each of these cases, it is important to develop methodologies that support the capability to design and operate microgrids in the absence of transmission over long periods of time. Unfortunately, such planning problems tend to be computationally difficult to solve and those that are straightforward to solve often lack the modeling fidelity that inspires confidence in the results. To address these issues, we first develop a high fidelity model for design and operations of a microgrid that include component efficiencies, component operating limits, battery modeling, unit commitment, capacity expansion, and power flow physics; the resulting model is a mixed-integer quadratically-constrained quadratic program (MIQCQP). We then develop an iterative algorithm, referred to as the Model Predictive Control (MPC) algorithm, that allows us to solve the resulting MIQCQP. We show, through extensive computational experiments, that the MPC-based method can scale to problems that have a very long planning horizon and provide high quality solutions that lie within 5\% of optimal.

\end{abstract}



\section{Introduction}
\label{Sec:intro}

In recent years, there has been increasing concerns about the susceptibility of modern electric power systems to extreme events that cause large-scale black outs. Even in 2017, hurricanes such as Harvey, Irma, and Maria have demonstrated the ability of such events to leave populations in the United States and the Caribbean without power for days or even weeks. As a result, the power engineering community has a vested interest in exploring solutions that can mitigate the impacts of these events. One proposed solution is microgrids, i.e. a decentralized subsystem that is capable of operating in either a grid-connected mode or an islanded mode to satisfy the power demand of a local community during large-scale black outs.

One of the biggest obstacles for rapid deployment of microgrids is their high initial outlay. 
This makes optimizing the design and operation of microgrids an economically very important problem. 
However, this problem is generally very difficult to solve because of the underlying nonlinear, nonconvex physics and discrete decision variables. The challenges in solving this problem are further exacerbated by desires to ensure that islanded microgrids can operate for long periods of time without a connection to bulk transmission systems.
In spite of these difficulties, recent papers \cite{madathil2017resilient,mashayekh2017security,barnes2017tools,DERCAM2015,byeon2018communication} have considered and modeled various important details of this problem, \textit{i.e.} communication-constrained expansion planning, component efficiency, linearized AC physics, and N-1 security. However, the problem of planning and operating a microgrid over a finite time horizon, \textit{i.e.} months, or even years, \cite{ton2012us}, remains intractable. Even planning for peak or average load conditions (\textit{i.e.} typical-day) is computationally very difficult.
 
In this paper, we focus on the computational tractability of long-time horizon planning and design of microgrids and push the boundaries of scalability beyond what is possible for existing approaches. In this context, we remark that for the purposes of this article, we have ignored uncontrollable sources like wind and solar in order to isolate and examine the difficulty in developing efficient methods. 
Traditionally, long-time horizon problems in the control literature are handled by decomposing the time-horizon into smaller time-stages, i.e. receding-horizon (RH) methods \cite{camacho2013model,scattolini2007hierarchical}.
In RH, each sub-problem is sequentially solved to optimality and the solutions from previous stages serve as initial conditions for subsequent stages. Generally speaking, this approach is computationally attractive and such methods have been applied to microgrid design \cite{bidram2012hierarchical,parisio2014model,fortenbacher2017optimal}.
However, optimality of individual stages does not guarantee global optimality, even when the problem is convex \cite{boyd2011distributed}. Furthermore, the RH approach is short-sighted in the sense that it does not take into account for any fluctuations in the load (or the slow dynamics) in the subsequent stages.  To address this issue, recent work has developed a first-order method that uses dual information from a coarse long-horizon problem \cite{zavala2016new} to communicate information between sub problems. Under certain assumptions, such as convexity, this method is guaranteed to converge to global optimality and has been shown to be effective on long-time horizon, continuous problems.

In contrast, high fidelity models of microgrid design that include discrete installation of technologies, AC power flow physics, engineering limits, and battery storage device efficiencies is a long-time horizon mixed-integer, nonlinear problem (MINLP) \cite{madathil2017resilient}. To address these key modeling details, this paper generalizes the existing first-order methods \cite{zavala2016new} to a mixed-integer setting. In short, the key contributions of this paper are (i) formulation of the Microgrid Design and Operation Problem ({\problem}) as a mixed-integer quadratically constrained quadratic program (MIQCQP), (ii) a convex relaxation of battery-efficiency modeling without the use of disjunctive binary variables, (iii) an iterative first-order-based MPC algorithm to solve to the {\problem } for long-time horizon; the algorithm is very general and can be applied to any MIQCQP that permits a decomposition in space or time, (iv) an extensive computational experimentation of the algorithm on an IEEE 13-node test feeder to corroborate its effectiveness in finding high quality solutions (typically, within 5\% of the optimal) for the \problem. 

The rest of the paper is organized as follows: the next section details the nomenclature and presents the MIQCQP model for the \problem. Sec. \ref{sec:algo} presents the MPC algorithm by starting with a QCQP and then extending it to an MIQCQP. Finally, Sec. \ref{sec:results} presents extensive computational results followed by conclusions and future work in Sec. \ref{sec:conc}.

\section{Mathematical formulation \label{sec:formulation}}

\subsection{Nomenclature} \label{subsec:nomenclature}
Throughout the rest of the article, boldface symbols are used to denote constants. \\
\noindent \textbf{Sets}: \\
$\mathcal N$ - set of nodes (buses), indexed by $i$  \\
$\mathcal E$ - set of \textit{from} edges (lines), indexed by $(i,j)$\\
$\mathcal E^r$ - set of \textit{to} edges (lines) \\
$\mathcal T$ - set of time periods, indexed by $t$ \\
$\mathcal B$ - set of batteries, indexed by $b$ \\
$\mathcal D$ - set of diesel generators, indexed by $d$ \\
$\mathcal D(i)$ - set of diesel generators at bus $i$\\
$\mathcal B(i)$ - set of batteries at bus $i$ \\
$\mathcal R = \mathcal B \cup \mathcal D$ - set or resources, indexed by $r$ \\
$\mathcal R(i)$ - set of resources at bus $i$\\

\noindent \textbf{Constants}: \\
$\bm f_b,\bm f_d$ - fixed cost for $b \in \mathcal B$, $d \in \mathcal D$, \$ \\
$\bm g_b$ - variable capacity cost for $b \in \mathcal B$, \$/MVA \\
$\bm c^0_d, \bm c_d^1, \bm c_d^2$ - operational cost coefficients for $d \in \mathcal D$\\
$\bm{UT}_d, \bm{DT}_d$ - minimum up-time and down-time for $d \in \mathcal D$\\
$\bm{RU}_d, \bm{RD}_d$ - ramp-up and ramp-down limit for $d \in \mathcal D$\\
$\bm s_{ij}$ - thermal limit for line $(i,j)$, MVA \\
$\bm{lp}_i^t + \bm j \bm{lq}_i^t$ - AC power demand at bus $i \in \mathcal N$, MVA \\
$\bm r_{ij} + \bm j \bm x_{ij}$ - impedance of line $(i,j) \in \mathcal E$ \\
$\bm p_{r}^{gl}, \bm p_{r}^{gu}$ - real power generation limits for $r \in \mathcal R$, MW \\
$\bm q_r^{gl}, \bm q_{r}^{gu}$ - reactive power generation limits for $r \in \mathcal R$, MVAr \\
$\bm v_i^l, \bm v_i^u$ - bounds for voltage magnitude square at bus $i$, kV \\
$\bm m_b$ - maximum installable capacity for $b\in \mathcal B$, MVA \\
$\bm{\overline{sc}}_b$ - maximum energy storage capacity of $b \in \mathcal B$, MW-hr \\
$\bm{\uptau}_b$ - energy storage capacity of $b \in \mathcal B$, MVA \\ 
$\bm{\eta}_{\text{ch}}^b, \bm{\eta}_{\text{dis}}^b$ - charging, discharging efficiency for battery $b \in \mathcal B$ \\
$\bm{\eta}^d$ - efficiency of diesel generator $d \in \mathcal D$\\
$\bm \mu$ - load shedding penalty, \$/MW \\
$\bm h_b(i)$ - number of batteries that can be built at bus $i$ \\
$\bm h_d(i)$ - number of generators that can be built at bus $i$ \\
$\bm \Delta \bm t$ - duration of time step, hrs \\

\noindent \textbf{Build decision variables}: \\
$z_b$ - binary build variable for each battery $b \in \mathcal B$ \\
$z_d$ - binary build variable for each generator $d \in \mathcal D$ \\

\noindent \textbf{Operation variables for diesel generators}: \\
$x_d^t$ - active/inactive status of $d \in D$ for $t \in \mathcal T$\\
$y_d^t$ - start-up status of $d \in \mathcal D$ for $t \in \mathcal T$\\
$w_d^t$ - shut-down status of $d \in \mathcal D$ for $t \in \mathcal T$\\ 
$p_d^{gt} + \bm j q_d^{gt}$ - apparent power generated by $d \in \mathcal D$ for $t \in \mathcal T$\\
$\widehat{p}_d^{gt}$ - active power generated by $d \in \mathcal D$ for $t\in \mathcal T$ before losses\\

\noindent \textbf{Operation variables for batteries}: \\
$p_b^{gt} + \bm j q_b^{gt}$ - apparent power entering/leaving $b \in \mathcal B$ for $t \in \mathcal T$\\
$\widehat{p}_b^{gt}$ - active power stored in $b \in \mathcal B$ for $t\in \mathcal T$\\
$sc_b^t$ - state-of-charge (energy stored) in $b \in \mathcal B$ at $t\in \mathcal T$, MW-hr \\
$s_b$ - maximum apparent power generation for $b \in \mathcal{B}$, MVA\\

\noindent \textbf{Other variables}: \\
$v_i^t$ - squared voltage on node $i \in \mathcal N$ at $t \in \mathcal T$ \\
$p_e^t, q_e^t$ - active/reactive power flow on line $e \in \mathcal E$ at $t\in \mathcal T$ \\
$\ell_i^{pt}, \ell_i^{qt}$ - active/reactive load shed on node $i \in \mathcal N$ at $t \in \mathcal T$ 

\subsection{Optimization problem} \label{subsec:opt}
With the notations presented in Sec. \ref{subsec:nomenclature}, we next formulate a mixed-integer quadratically constrained quadratic programming formulation for the \problem. The objective function of the {\problem }, Eq. \eqref{objective}, minimizes the build cost, the total load shed, and the operating cost of the resources (batteries and diesel generators). The operating cost of diesel generators includes a no-load cost and a quadratic generation cost.
\begin{equation}
\begin{split}
\min  \, & \sum_{b \in \mathcal B} \left( \bm f_b z_b + \bm g_b {s}_b \right) + \sum_{d \in \mathcal D} \bm f_d z_d + \sum_{t \in \mathcal T} \sum_{i \in \mathcal N} \bm \mu \left(\ell_i^{pt} + \ell_i^{qt}\right) \\
&  + \sum_{t \in \mathcal T} \sum_{d \in \mathcal D} \left[ \bm c_d^0 \cdot x_d^t + \bm c_d^1 \cdot \widehat{p}_d^{gt} +  \bm c_d^2 \cdot \left(\widehat{p}_d^{gt} \right)^2\right]
\end{split} 
\label{objective}
\end{equation}
The {\problem } is subject to the following constraints: 

\noindent \textit{\underline{\smash{\textbf{Power flow physics}}}} \\ The physics of power flow, thermal limits, and voltage bounds are described by the following set of constraints:
\begin{subequations}
\begin{flalign}
& \sum_{r \in \mathcal R(i)} p_r^{gt}+ \ell_i^{pt} -\bm{lp}_i^t = \sum_{(i,j) \in \mathcal E\cup \mathcal E^r} p_{ij}^t \quad \forall i\in \mathcal N, t \in \mathcal T, & \label{eq:powerflow1} \\
&  \sum_{r \in \mathcal R(i)} q_r^{gt} + \ell_i^{qt} -\bm{lq}_i^t = \sum_{(i,j) \in \mathcal E\cup \mathcal E^r} q_{ij}^t \quad \forall i\in \mathcal N, t \in \mathcal T, & \label{eq:powerflow2} \\
& v_j^t = v_i^t - 2 \cdot (\bm r_{ij} p_{ij}^t + \bm x_{ij} q_{ij}^t) \quad \forall (i,j) \in \mathcal E, t\in \mathcal T, & \label{eq:powerflow3} \\
& \left( p_{ij}^t \right)^2 + \left( q_{ij}^t \right)^2 \leqslant \bm s_{ij}^2 \quad \forall (i,j) \in \mathcal E, t\in \mathcal T, & \label{eq:thermallimit} \\
& \bm v_i^l \leqslant v_i^t \leqslant  \bm v_i^u \quad \forall i \in \mathcal N, t \in \mathcal T. & \label{eq:voltagebound}
\end{flalign}
\end{subequations}
Eqs. \eqref{eq:powerflow1}--\eqref{eq:powerflow2} enforce Kirchoff's current law and Eq. \eqref{eq:powerflow3} models the ``LinDistFlow'' approximation of the nonconvex DistFlow AC power flow equations \cite{gan2015exact}. 
Eqs. \eqref{eq:thermallimit} and \eqref{eq:voltagebound} enforce the thermal limit for each transmission line and bounds on the voltage magnitude at each bus, respectively. 

\noindent \textit{\underline{\smash{\textbf{Resource limits}}}} \\ The following constraints enforce the capacity limits on the build and generation decisions for all the resources in $\mathcal R$.
\begin{subequations}
\begin{flalign}
& \sum_{b \in \mathcal B(i)} z_b \leqslant \bm h_b(i) \quad \forall i \in \mathcal N, & \label{eq:build_b} \\
& \sum_{d \in \mathcal D(i)} z_d \leqslant \bm h_d(i) \quad \forall i \in \mathcal N, & \label{eq:build_d} \\
& s_b \leqslant z_b \bm m_b \quad \forall b \in \mathcal B. & \label{eq:capacitylimit_battery}
\end{flalign}
\end{subequations}
Eqs. \eqref{eq:build_b} and \eqref{eq:build_d} restrict the number of batteries and diesel generators that can be installed at each bus, respectively. Eq. \eqref{eq:capacitylimit_battery} limits the capacity of each battery.

\noindent \textit{\underline{\smash{\textbf{Operation of diesel generators}}}} \\ 
We first introduce sets
$\mathcal A(d,t)$ and $\mathcal B(d,t)$ that support the formulation of minimum up-time and down-time constraints for each diesel generator $d\in D$ at each time period $t \in \mathcal T$.
\begin{align*}
& \mathcal A(d,t) = \{ \tilde{t} \in \mathcal T : t - \bm{UT}_d + 1 \leqslant \tilde t \} \,\, \forall d\in \mathcal D, t \in \mathcal T, & \\
& \mathcal B(d,t) = \{ \tilde{t} \in \mathcal T : t - \bm{DT}_d + 1 \leqslant \tilde t \} \,\, \forall d\in \mathcal D, t \in \mathcal T. & 
\end{align*}
The constraints on diesel generator operation are then:
\begin{subequations}
\begin{flalign}
&x_d^t \leqslant z_d \quad \forall d \in \mathcal D, t \in \mathcal T, & \label{eq:logic1} \\
&y_d^t - w_d^t = x_d^t - x_d^{t-1} \quad \forall d \in \mathcal D, t\in \mathcal T, & \label{eq:logic2}\\
&y_d^t + w_d^t \leqslant 1 \quad \forall d \in \mathcal D, t \in \mathcal T, & \label{eq:logic3} \\
& x_d^t \,\bm p_d^{gl} \leqslant \widehat{p}_d^{gt} \leqslant x_d^t \,\bm p_d^{gu} \quad \forall d \in \mathcal D, t\in \mathcal T, &\label{eq:genlimitp_d} \\
& x_d^t \,\bm q_d^{gl} \leqslant q_d^{gt} \leqslant x_d^t \,\bm q_d^{gu} \quad \forall d \in \mathcal D, t\in \mathcal T, &\label{eq:genlimitq_d} \\
& p_d^{gt} = \bm \eta^d \widehat{p}_d^{gt} \quad \forall d \in \mathcal D, t \in \mathcal T, & \label{eq:efficiency_d} \\
&\bm{RD}_d \geqslant p_d^{g(t-1)} - p_d^{gt} \quad \forall d \in \mathcal D, t\in \mathcal T, & \label{eq:rampdown}\\
&\bm{RU}_d \geqslant p_d^{gt} - p_d^{g(t-1)} \quad \forall d \in \mathcal D, t\in \mathcal T, & \label{eq:rampup}\\
&\sum_{k \in \mathcal A(d,t)} y_{d}^k \leqslant x_{d}^t \quad \forall d \in \mathcal D, t \in \mathcal T, \label{eq:uptime} \\
&\sum_{k \in \mathcal B(d,t)} w_{d}^k \leqslant 1-x_{d}^t \quad \forall d \in \mathcal D, t \in \mathcal T. \label{eq:downtime}
\end{flalign}
\end{subequations}
Eq. \eqref{eq:logic1} forces a diesel generator to be committed only if it is built. Eq. \eqref{eq:logic2} determines if the generator is started up or shut down at time period $t$ based of its on-off status between time periods $t$ and $t-1$. Eq. \eqref{eq:logic3} ensures that a generator $d$ is not started up and shut down in the same time period, $t$. Eqs. \eqref{eq:genlimitp_d} and \eqref{eq:genlimitq_d} enforce the generation limits of diesel generators. 
Real power limits are applied to $\widehat{p}_d^{gt}$, i.e. the generation produced before losses.
This loss is quantified by Eq. \eqref{eq:efficiency_d}. Eqs. \eqref{eq:rampdown} and \eqref{eq:rampup} enforce the ramping limits on consecutive time periods on every generator $d \in \mathcal D$. Finally, Eqs. \eqref{eq:uptime} and \eqref{eq:uptime} enforce the minimum up time and minimum down time constraints for every diesel generator. Usage of distinct binary variables, ($x^t_d,y^t_d,z^t_d$), to represent a tighter unit commitment polytope is akin to the formulation described in \cite{sundar2016unit}.  

\noindent \textit{\underline{\smash{\textbf{Battery operation}}}}\\ The constraints 
on battery operations are given by:
\begin{subequations}
\begin{flalign}
& (p_b^{gt})^2 + (q_b^{gt})^2 \leqslant (s_b)^2 \quad \forall b\in \mathcal B, t \in \mathcal T &  \label{eq:battery_limit} \\
& sc_b^t = sc_b^{t-1} - \widehat{p}_b^{gt} \cdot \bm \Delta \bm t, \quad \forall b \in \mathcal B, t \in \mathcal T, & \label{eq:state_of_charge} \\
& 0 \leqslant sc_b^t \leqslant z_b \bm{\overline{sc}}_b \quad \forall b\in \mathcal B, t \in \mathcal T. \label{eq:capacity_limits} 
\end{flalign}
\end{subequations}
Eq. \eqref{eq:battery_limit} enforces the apparent power limits on charging and discharging for each battery. Eq. \eqref{eq:state_of_charge} determines the current state-of-charge of the battery, $sc_b^t$, based on whether the battery is being charged ($\widehat{p}_b^{gt} < 0$) or discharged ($\widehat{p}_b^{gt} > 0$). Finally, Eq. \eqref{eq:capacity_limits} imposes bounds on the energy storage of the battery. The model of charging and discharging for batteries is adapted from \cite{koutsopoulos2011optimal}.

\begin{figure}[!h]
\centering
\includegraphics[scale=0.33]{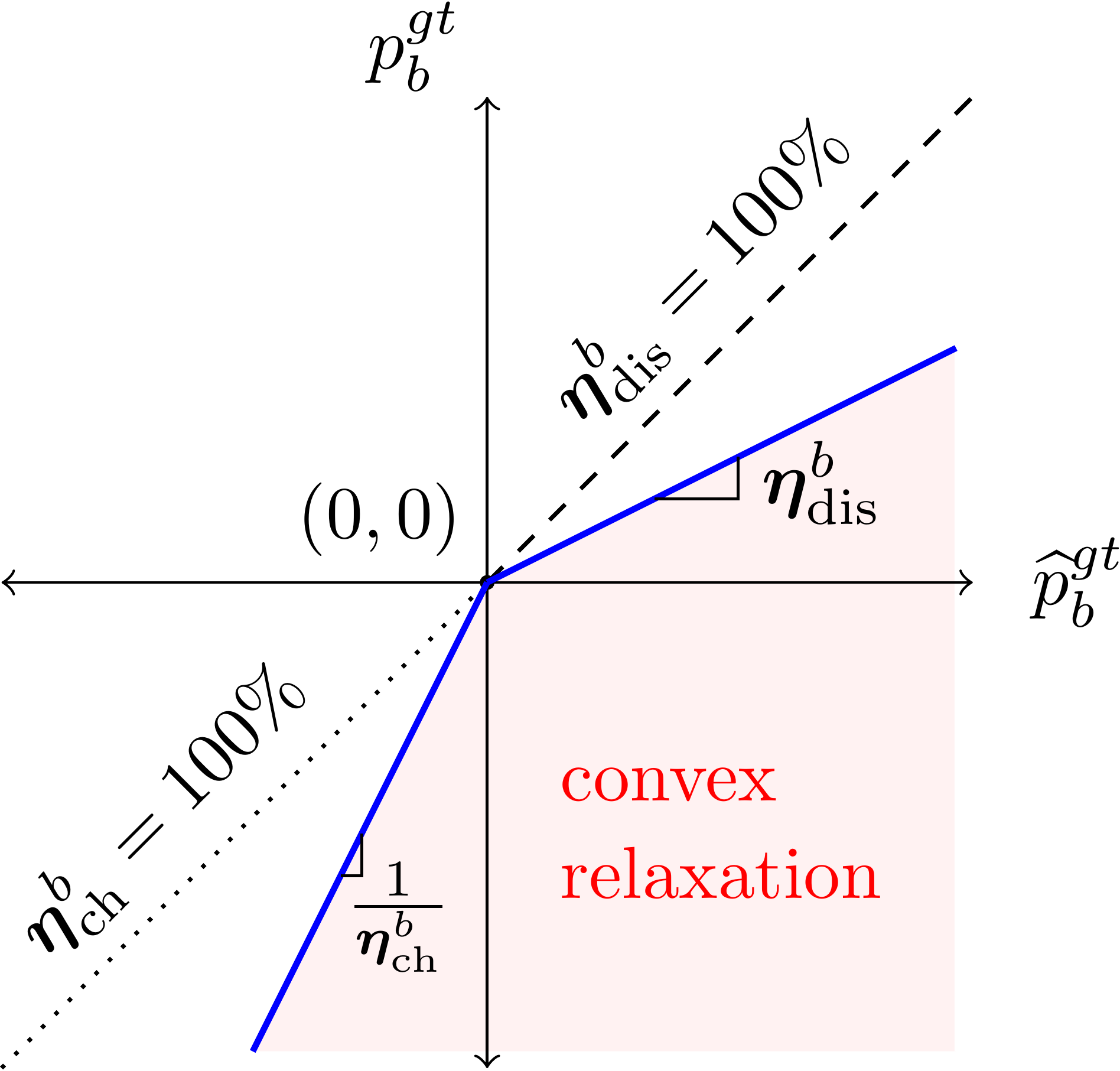}
\caption{An illustrative example of convex relaxation (shaded) of piecewise linear efficiency curve (blue) for batteries.}
    \label{fig:eff_modeling}
\end{figure}

\noindent \textit{\underline{\smash{\textbf{Battery efficiency modeling}}}} \\ We next discuss the constraints that model battery's charging and discharging efficiency without the use of disjunctive binary variables, typically used in the literature \cite{}. A preliminary version of the battery model, not utilizing disjunctive binary variables, without a formal explanation has been developed in \cite{madathil2017resilient}. In this article, we formalize the model as a convex relaxation. For ease of exposition, we present battery efficiency modeling with one charging and one discharging efficiency values, $\bm \eta_{\operatorname{ch}}^b$ and $\bm \eta_{\operatorname{dis}}^b$, respectively. Without loss of generality, this model is extendable to piece-wise linear charging and discharging curves with monotonically decreasing charging and discharging efficiency values (see \cite{madathil2017resilient}). We first define the charging and discharging regimes for battery $b \in \mathcal B$ and time period $t \in \mathcal T$ as:
\begin{subequations}
\begin{flalign}
& \text{discharging: } p_b^{gt} > 0, \, \widehat{p}_b^{gt} > 0,\, |\widehat{p}_b^{gt}| > |p_b^{gt}|, & \label{eq:discharging_regime} \\
& \text{charging: } p_b^{gt} < 0, \,\widehat{p}_b^{gt} < 0,\, |\widehat{p}_b^{gt}| < |p_b^{gt}|. & \label{eq:charging_regime} 
\end{flalign}
\end{subequations}
For any battery $b \in \mathcal B$ and time period $t \in \mathcal T$, $\widehat{p}_b^{gt}$ is the power stored in the battery and $p_b^{gt}$ is the amount of power entering (charging) or leaving (discharging) $b$. Then, the battery efficiency is modeled using this equation:
\begin{flalign}
& \begin{pmatrix} 
p_b^{gt} = \bm \eta_{\operatorname{dis}}^b\, \widehat{p}_b^{gt} \\
p_b^{gt},  \,\widehat{p}_b^{gt} > 0, \\
|\widehat{p}_b^{gt}| > |p_b^{gt}| 
\end{pmatrix}  \lor 
\begin{pmatrix} p_b^{gt} = \frac 1 {\bm \eta_{\operatorname{ch}}^b} \,\widehat{p}_b^{gt} \\
p_b^{gt}, \,\widehat{p}_b^{gt} < 0, \\
|\widehat{p}_b^{gt}| < |p_b^{gt}|
\end{pmatrix}
\quad \forall b\in \mathcal B, t \in \mathcal T. &\label{eq:eff_and} 
\end{flalign}
Traditionally, Eq. \eqref{eq:eff_and} is reformulated into disjunctive, linear constraints (i.e., \cite{felder2014optimal}). In this article, we introduce the following convex relaxation: 
\begin{subequations}
\begin{flalign}
& p_b^{gt} \leqslant \bm \eta_{\operatorname{dis}}^b \,\widehat{p}_b^{gt} \quad \forall b \in \mathcal B, t \in \mathcal T, & \label{eq:discharging} \\
& p_b^{gt} \leqslant \frac 1 {\bm \eta_{\operatorname{ch}}^b} \,\widehat{p}_b^{gt}  \quad \forall b \in \mathcal B, t \in \mathcal T. & \label{eq:charging}
\end{flalign}
\end{subequations}
The convex relaxation of Eq. \eqref{eq:eff_and} is given by Eqs. \eqref{eq:discharging}, \eqref{eq:charging}, \eqref{eq:battery_limit}, and \eqref{eq:capacity_limits}. This convex relaxation of the efficiency curves is also illustrated in Fig. \ref{fig:eff_modeling}. Though we assume the piecewise efficiency curve with only two pieces, generalization of the proposed convex relaxation approach is straightforward for the case when there are multiple charging and discharging efficiencies (see \cite{madathil2017resilient}). We remark that although Eqs. \eqref{eq:discharging} and \eqref{eq:charging} yield a relaxation, the relaxed solutions are often empirically tight (see the computational results).

\section{Algorithms}
\label{sec:algo}
The formulation presented for the {\problem } is an MIQCQP. A typical planning horizon for the {\problem } ranges from days to a few years. For the purposes of this paper, we assume the variations in load over a full year can be reduced to a few weeks, based on the typical-day approach \cite{madathil2017resilient,mashayekh2017security},
with a time discretization of $15$ minutes. Using state-of-the-art solvers such as Gurobi/CPLEX, the {\problem }, even for a time horizon of $1$ day at $15$ minute time intervals, is challenging to solve to optimality 
because of the integer variables.  
To address the computational difficulty of the \problem, this section describes an efficient heuristic that adapts hierarchical model-predictive control architectures \cite{zavala2016new} to discrete problems. 
We first present the MPC algorithm for a continuous relaxation of the \problem, i.e. a QCQP, and then extend it to MIQCQPs. For ease of exposition, we recast the formulation of Sec. \ref{subsec:opt} as:
\begin{subequations}
\begin{flalign}
& (\mathcal F_1): \quad \min \sum_{t\in \mathcal T} \varphi (x_t, y_t) \text{ subject to:}& \label{eq:obj} \\
&  g(x_{\tau}, y_{\tau}) = 0 \quad \text{where } \tau \subseteq \{0, \dots, |\mathcal T|\}, & \label{eq:cons1} \\
& x_0 = \overline{\bm x}, y_0 = \overline{\bm y}, & \label{eq:cons2} \\
& x_t \in \mathbb R^n, y_t \in \{0,1\}^m \quad \forall t \in \mathcal T. & \label{eq:var_bounds}
\end{flalign}
\end{subequations}
The variables $x_t$ and $y_t$ model the vectors of continuous and binary variables in {\problem}'s formulation. Eq. \eqref{eq:cons1} defines the feasible region of the optimization problem, as specified by the linear and quadratic constraints in Sec. \ref{subsec:opt}. Eq. \eqref{eq:cons2} specifies the initial condition and Eq. \eqref{eq:var_bounds} specifies the variable domains. 

One of the main challenges with solving {\problem } to optimality is the large cardinality of $\mathcal T$. To address this challenge, we partition the time horizon $\{0,1,.., |\mathcal T|\}$ into $S$ stages. Each stage consists of a subset of consecutive time periods. For ease of exposition, we assume, without loss of generality, that every stage has $K$ time periods numbered, $\mathcal K = \{1,.., K\}$. By convention, the zeroth time period for every stage defines the initial conditions for that stage. We let $\mathcal S = \{1,.., S\}$. Using this notation, $\mathcal F_1$ is reformulated as $\mathcal F_2$:
\begin{subequations}
\begin{flalign}
& (\mathcal F_2): \quad \min \sum_{s \in \mathcal S} \sum_{k \in \mathcal K} \varphi (x_{s,k}, y_{s,k}) \text{ subject to:} & \label{eq:obj_1} \\
& g(x_{s,\tau}, y_{s,\tau}) = 0 \quad \forall s \in \mathcal S, \tau = \mathcal K,  & \label{eq:cons1_1}\\
& x_{0,K} = \overline{\bm x}, y_{0,K} = \overline{\bm y}, & \label{eq:cons2a_1} \\
& x_{s,0} = x_{s-1,K}, \, y_{s,0} = y_{s-1,K} \quad \forall s \in \mathcal S, k \in \mathcal K, & \label{eq:cons2b_1}\\
& x_{s,k} \in \mathbb R^n, y_{s,k} \in \{0,1\}^m \quad \forall s \in \mathcal S, k \in \mathcal K. & \label{eq:var_bounds_1}
\end{flalign}
\end{subequations}
In formulation $\mathcal F_2$, Eq. \eqref{eq:cons1_1} defines all the constraints of stage $s \in \mathcal S$. These constraints are referred to as ``stage constraints''. Eqs. \eqref{eq:cons2a_1} and \eqref{eq:cons2b_1} specify the initial conditions for each stage. Specifically, Eq. \eqref{eq:cons2a_1} specifies the initial conditions for the first stage and Eq. \eqref{eq:cons2b_1} states that the initial condition of stage $s \in \mathcal S$ is equal to the $K$\textsuperscript{th} (final) time period of the previous stage. The constraints in Eq. \eqref{eq:cons2b_1} are referred to as ``coupling constraints''. Finally, Eq. \eqref{eq:var_bounds_1} specifies the binary restrictions on the variables. 

\subsection{MPC algorithm for the continuous relaxation of the \problem} 
\label{subsec:cont_relaxation}
To simplify the presentation of the algorithm, without loss of generality, we ignore the stage constraints and focus on the coupling constraints. We first group the variables by stages by defining vectors of variables $\vec{x}_s = (x_{s,0}, x_{s,1}, \dots, x_{s,K})$ and $\vec y_s = (y_{s,0}, y_{s,1}, \dots, y_{s,K})$. The block form of $\mathcal F_2$, without the stage constraints, is then defined as:
\begin{subequations}
\begin{flalign}
& (\mathcal F_3): \quad \min \sum_{s \in \mathcal S} \left( \frac 12 \begin{bmatrix}  \vec x_s \\ \vec y_s \end{bmatrix}^\intercal \bm Q_s \begin{bmatrix}  \vec x_s \\ \vec y_s \end{bmatrix} - \bm c_s^{\intercal} \begin{bmatrix}  \vec x_s \\ \vec y_s \end{bmatrix} \right) & \label{eq:obj_2} \\
& (\lambda_s) \quad \overline{\bm \Pi}_s \begin{bmatrix} \vec x_s \\ \vec y_s \end{bmatrix} = \munderbar{\bm \Pi}_s \begin{bmatrix} \vec x_{s-1} \\ \vec y_{s-1} \end{bmatrix}, \quad \forall s \in \mathcal S, \label{eq:coupling} & \\
& \vec x_s \in \mathbb R^{n \times K}, \vec y_s \in \{0,1\}^{m \times K} \quad \forall s \in \mathcal S. \label{eq:bin} &
\end{flalign}
\end{subequations}
Eq. \eqref{eq:coupling} is equivalent to Eqs. \eqref{eq:cons2a_1} and \eqref{eq:cons2b_1} in matrix form. The $\lambda_s$ are the dual variables of  constraints \eqref{eq:coupling}. The matrices $\munderbar{\bm \Pi}_s$ and $\overline{\bm \Pi}_s$ are the coefficient matrices; they correspond to the constraints \eqref{eq:cons2b_1}, rewritten for each stage $s \in \mathcal S$ in a matrix form. The variables for $s = 0$ are fixed to the initial conditions of the \problem. When the binary constraints of
 $\mathcal F_3$ are relaxed, the first-order optimality conditions (KKT) are given by
\begin{subequations}
\begin{flalign}
& \bm Q_s \begin{bmatrix} \vec x_s \\ \vec y_s \end{bmatrix} - \overline{\bm \Pi}_s^{\intercal} \lambda_s + \munderbar{\bm \Pi}_{s+1}^{\intercal} \lambda_{s+1} = \bm c_s \quad \forall s \in \mathcal S\setminus \{S\}, & \label{eq:kkt1} \\
& \bm Q_s \begin{bmatrix} \vec x_s \\ \vec y_s \end{bmatrix} - \overline{\bm \Pi}_s^{\intercal} \lambda_s  = \bm c_s \quad \forall s \in \{S\}, & \label{eq:kkt2} \\
& \overline{\bm \Pi}_s \begin{bmatrix} \vec x_s \\ \vec y_s \end{bmatrix} = \munderbar{\bm \Pi}_s \begin{bmatrix} \vec x_{s-1} \\ \vec y_{s-1} \end{bmatrix} \quad \forall s \in \mathcal S. \label{eq:kkt3}
\end{flalign}
\label{eq:kkt}
\end{subequations}
Any technique to solve the KKT system given in Eq. \eqref{eq:kkt} is iterative (for instance, Gauss Seidal). 
To avoid centrally solving this KKT system using an iterative method, we use a decentralized technique with an update index $\ell \in \mathbb Z_+$; the update equations used for each iteration is given below:
\begin{subequations}
\begin{flalign}
& \begin{bmatrix} 
\bm Q_1 & - \overline{\bm \Pi}_1^{\intercal} \\
\munderbar{\bm \Pi}_{1} & \bm 0
\end{bmatrix}
\begin{bmatrix}
\vec x_1^{\ell+1} \\ \vec y_1^{\ell+1} \\ \lambda_1^{\ell+1} 
\end{bmatrix} = 
\begin{bmatrix}
\bm c_1 \\ 0
\end{bmatrix} -
\begin{bmatrix}
\bm 0 & \munderbar{\bm \Pi}^{\intercal}_{2} \\
\overline{\bm \Pi}_{1} & \bm 0
\end{bmatrix}
\begin{bmatrix}
\vec x_{0}^{\ell} \\ \vec y_{0}^{\ell} \\ \lambda_{2}^{\ell} 
\end{bmatrix} & \label{eq:dec0} \\
& \begin{bmatrix} 
\bm Q_s & - \overline{\bm \Pi}_s^{\intercal} \\
\munderbar{\bm \Pi}_{s} & \bm 0
\end{bmatrix}
\begin{bmatrix}
\vec x_s^{\ell+1} \\ \vec y_s^{\ell+1} \\ \lambda_s^{\ell+1} 
\end{bmatrix} = 
\begin{bmatrix}
\bm c_s \\ 0
\end{bmatrix} -
\begin{bmatrix}
\bm 0 & \munderbar{\bm \Pi}^{\intercal}_{s+1} \\
\overline{\bm \Pi}_{s} & \bm 0
\end{bmatrix}
\begin{bmatrix}
\vec x_{s-1}^{\ell+1} \\ \vec y_{s-1}^{\ell+1} \\ \lambda_{s+1}^{\ell+1}
\end{bmatrix} & \notag \\ & \qquad \qquad \qquad \qquad \qquad \qquad \qquad \qquad \forall s \in \mathcal S\setminus \{1,S\}, & \label{eq:dec1} \\
& \begin{bmatrix} 
\bm Q_S & - \overline{\bm \Pi}_S^{\intercal} \\
\munderbar{\bm \Pi}_{S} & \bm 0
\end{bmatrix}
\begin{bmatrix}
\vec x_S^{\ell+1} \\ \vec y_S^{\ell+1} \\ \lambda_S^{\ell+1} 
\end{bmatrix} = 
\begin{bmatrix}
\bm c_S \\ 0
\end{bmatrix} -
\begin{bmatrix}
\bm 0 & \bm 0 \\
\overline{\bm \Pi}_{S} & \bm 0
\end{bmatrix}
\begin{bmatrix}
\vec x^{\ell+1}_{s-1} \\ \vec y^{\ell+1}_{s-1} \\ \lambda^{\ell+1}_{S} 
\end{bmatrix}. & \label{eq:dec2} 
\end{flalign}
\label{eq:update_equations}
\end{subequations}
\begin{figure*}[h!]
    \centering
    \includegraphics[scale=0.9]{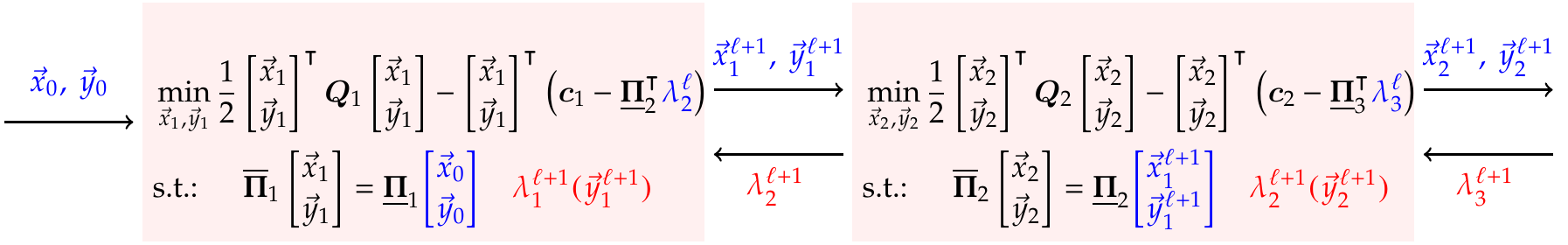}
    \caption{Sketch of the primal and dual variables being communicated to neighboring stages at iteration $(\ell+1)$.}
    \label{fig:pd_updates}
\end{figure*}
One key observation of Eq. \eqref{eq:update_equations} is that each of the decentralized iteration equations are first-order KKT conditions for the problems associated with each stage. The partitioned problems are defined as:
\begin{subequations}
\allowdisplaybreaks
\begin{flalign}
\begin{bmatrix} \vec x_1^{\ell+1} \\ \vec y_1^{\ell+1}\end{bmatrix} & = \operatornamewithlimits{argmin}_{\vec x_1, \vec y_1} 
 \frac 12 \begin{bmatrix}  \vec x_1 \\ \vec y_1 \end{bmatrix}^\intercal \bm Q_1 \begin{bmatrix}  \vec x_1 \\ \vec y_1 \end{bmatrix} -  \begin{bmatrix}  \vec x_1 \\ \vec y_1 \end{bmatrix}^\intercal \left(\bm c_1 - \munderbar{\bm \Pi}_2^{\intercal} \lambda_2^{\ell} \right)&  \notag \\
 & \text{subject to: } \quad  \overline{\bm \Pi}_1 \begin{bmatrix} \vec x_1 \\ \vec y_1 \end{bmatrix} = \munderbar{\bm \Pi}_1 \begin{bmatrix} \vec x_{0}^{\ell} \\ \vec y_{0}^{\ell} \end{bmatrix} \quad (\lambda_1^{\ell+1}). \label{eq:p1} &
\end{flalign}
\begin{flalign}
\begin{bmatrix} \vec x_s^{\ell+1} \\ \vec y_s^{\ell+1}\end{bmatrix} & =
\operatornamewithlimits{argmin}_{\vec x_s, \vec y_s} 
 \frac 12 \begin{bmatrix}  \vec x_s \\ \vec y_s \end{bmatrix}^\intercal \bm Q_s \begin{bmatrix}  \vec x_s \\ \vec y_s \end{bmatrix} -  \begin{bmatrix}  \vec x_s \\ \vec y_s \end{bmatrix}^\intercal \left(\bm c_s - \munderbar{\bm \Pi}_{s+1}^{\intercal} \lambda_{s+1}^{\ell} \right)& \notag \\
 & \text{subject to: } \quad \overline{\bm \Pi}_s \begin{bmatrix} \vec x_s \\ \vec y_s \end{bmatrix} = \munderbar{\bm \Pi}_s \begin{bmatrix} \vec x_{s-1}^{\ell+1} \\ \vec y_{s-1}^{\ell+1} \end{bmatrix} \quad (\lambda_s^{\ell+1}). \label{eq:ps} &
\end{flalign}
\begin{flalign}
\begin{bmatrix} \vec x_S^{\ell+1} \\ \vec y_S^{\ell+1}\end{bmatrix} & =
\operatornamewithlimits{argmin}_{\vec x_S, \vec y_S} 
 \frac 12 \begin{bmatrix}  \vec x_S \\ \vec y_S \end{bmatrix}^\intercal \bm Q_S \begin{bmatrix}  \vec x_S \\ \vec y_S \end{bmatrix} & \notag \\
 & \text{subject to: } \quad  \overline{\bm \Pi}_S \begin{bmatrix} \vec x_S \\ \vec y_S \end{bmatrix} = \munderbar{\bm \Pi}_S \begin{bmatrix} \vec x_{s-1}^{\ell+1} \\ \vec y_{s-1}^{\ell+1} \end{bmatrix} \quad (\lambda_S^{\ell+1}). \label{eq:pS}&
\end{flalign}
\label{eq:stage_problems}
\end{subequations}
Given $s \in \mathcal S$, we use $\mathcal P (\vec x_{s-1}^{\ell}, \vec y_{s-1}^{\ell}, \lambda_{s+1}^{\ell})$ to define the partitioned problem in Eq. \eqref{eq:stage_problems} corresponding to stage $s$.
Since the decentralized KKT system in \eqref{eq:update_equations} are first-order KKT conditions of $\mathcal P (\vec x_{s-1}^{\ell}, \vec y_{s-1}^{\ell}, \lambda_{s+1}^{\ell})$,
the KKT system is solved by solving the smaller (partitioned) optimization problems and communicating primal ($\vec x_s$, $\vec y_s$) and dual ($\lambda_s$) variables between partitions. Convexity of the relaxation of {\problem} (convex QCQP) guarantees convergence of this decentralized technique to the optimal solution (see \cite{zavala2016new}). 

\subsection{MPC algorithm for the \problem} \label{subsec:mpc_full}
The heuristic for solving the {\problem }, a MIQCQP, is derived by making the following two changes to the technique presented in Sec. \ref{subsec:cont_relaxation}: (i) the partitioned problem corresponding to stage $s \in \mathcal S$ in Eq. \eqref{eq:stage_problems} is solved as a MIQCQP and the primal solution ($\vec x_s$, $\vec y_s$) is communicated forward to the next stage (see step \ref{algo:1} of Algorithm \ref{algo:pseudocode} and Fig. \ref{fig:pd_updates}) and (ii) the dual values of the coupling constraints ($\lambda_s$) are computed by fixing the values of the binary variables, $\vec y_s$ to the solution of the MIQCQP for stage $s$ (see step \ref{algo:2} of Algorithm \ref{algo:pseudocode} and Fig. \ref{fig:pd_updates}). 
In other words, the dual value for iteration $\ell$ and stage $s$, $\lambda_s^{\ell+1}$, of the linking constraint is computed as a function of the binary variable solution values by solving the optimization problem $\mathcal P (\vec x_{s-1}^{\ell}, \vec y_{s-1}^{\ell}, \lambda_{s+1}^{\ell}) \cup \{\vec y_s = \vec y_s^{\ell+1} \}$. These values are denoted by $\lambda_s^{\ell+1}(\vec y_s^{\ell+1})$. Unlike solving the continuous relaxation of the {\problem}, this algorithm has no convergence guarantees to a global optimum. Hence, we stop the algorithm after a finite number of iterations. A pseudo-code of the MPC is shown in algorithm \ref{algo:pseudocode}. The update scheme at iteration $\ell$ is illustrated in the Fig. \ref{fig:pd_updates}. 

\begin{algorithm}
\caption{Psuedocode for the MPC algorithm}
\label{algo:pseudocode}
\small
\begin{algorithmic}[1]
\State $\ell \gets 1$
\State Solve continuous relaxation of {\problem } and initialize $\lambda_s^{\ell}$, for all $s \in \mathcal S$ \label{algo:3}
\While{$\ell \leqslant \bm N$}
    \For{$s=1, \dots, S$}
        \State $\vec x_{s}^{\ell+1}$, $\vec y_s^{\ell+1} \gets \mathcal P (\vec x_{s-1}^{\ell}, \vec y_{s-1}^{\ell}, \lambda_{s+1}^{\ell}) $ \label{algo:1}
        \State $\lambda_s^{\ell+1} \gets \operatornamewithlimits{getdual}\left(\mathcal P (\vec x_{s-1}^{\ell}, \vec y_{s-1}^{\ell}, \lambda_{s+1}^{\ell}) \bigcup \{\vec y_s = \vec y_s^{\ell+1} \}\right)$ \label{algo:2}
        \vspace{-0.3cm}
    \EndFor
    \State $\ell \gets \ell + 1$
\EndWhile
\end{algorithmic}
\end{algorithm}

\noindent \textit{\underline{\smash{Initialization of dual variable values:}}} 
Though there are many ways of initializing the dual values, setting them to zero values would be a trivial initialization. It is interesting to note that when the initial dual values are set to zero, the first iteration of MPC is equivalent to the standard \textit{Receding Horizon} (RH) algorithm \cite{zavala2016new}. However, one of the known drawbacks of RH is it's inability to capture the future time horizon's information in the current stage. To address this, MPC methods are applied in a more general setting of the \textit{hierarchical} multigrid control, which captures the slow-changing dynamics at low frequencies. In \problem, the real-time load fluctuations and slow raise to peak loads over 24 hours horizon can be viewed as the high and low frequency disturbances in the system, respectively. Thus, we initialize the duals of the coupling constraints by solving the continuous relaxation of {\problem} to optimality. This hierarchical MPC procedure (indirectly) captures and passes the low frequency fluctuations as cost-to-go values via the duals in the objective of each stage. Thus, MPC can be viewed as an iterative generalization of the RH scheme.



\section{Computational Results}
\label{sec:results}
We first present the test-system specifications and parameter values used in the formulation and in the MPC algorithm. 
\subsection{Data} \label{subsec:data}
\begin{figure}[h!]
    \centering
    \includegraphics[scale=0.56]{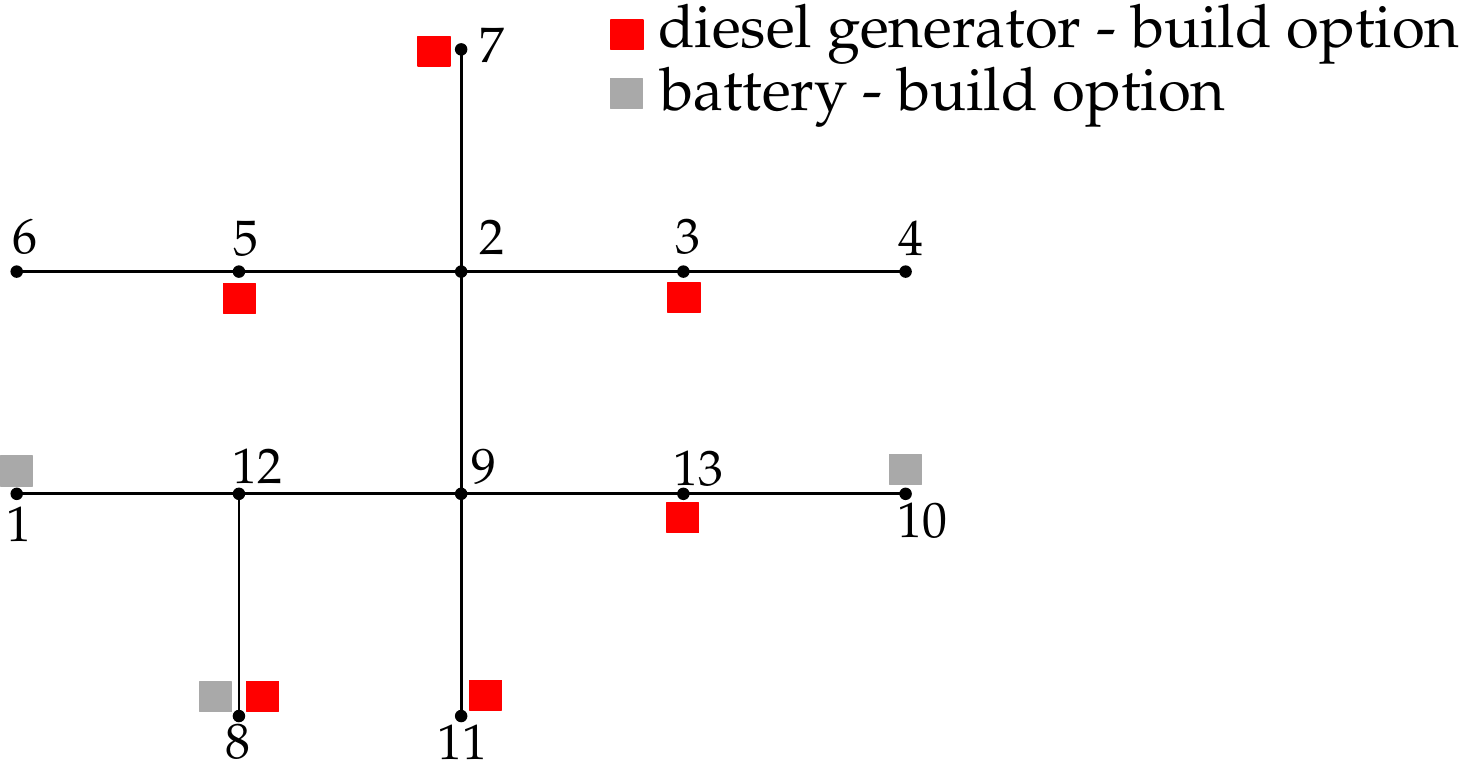}
    \caption{IEEE 13-node radial distribution network schematic.}
    \label{fig:schematic}
\end{figure}

\begin{figure}[h!]
  \centering
   \subfigure[Typical daily profile]{
   \includegraphics[scale=0.83]{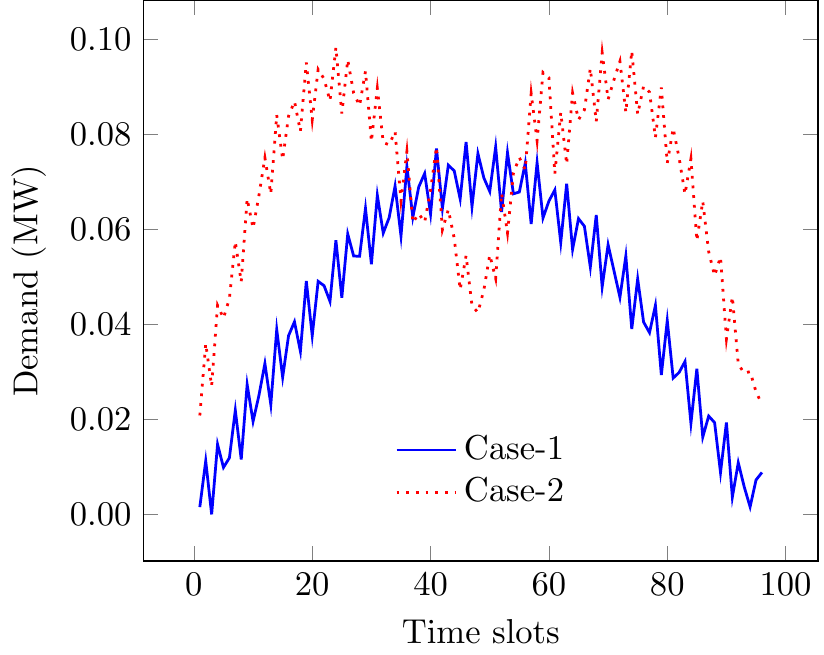}} 
   \subfigure[Typical weekly profile]{
   \includegraphics[scale=0.83]{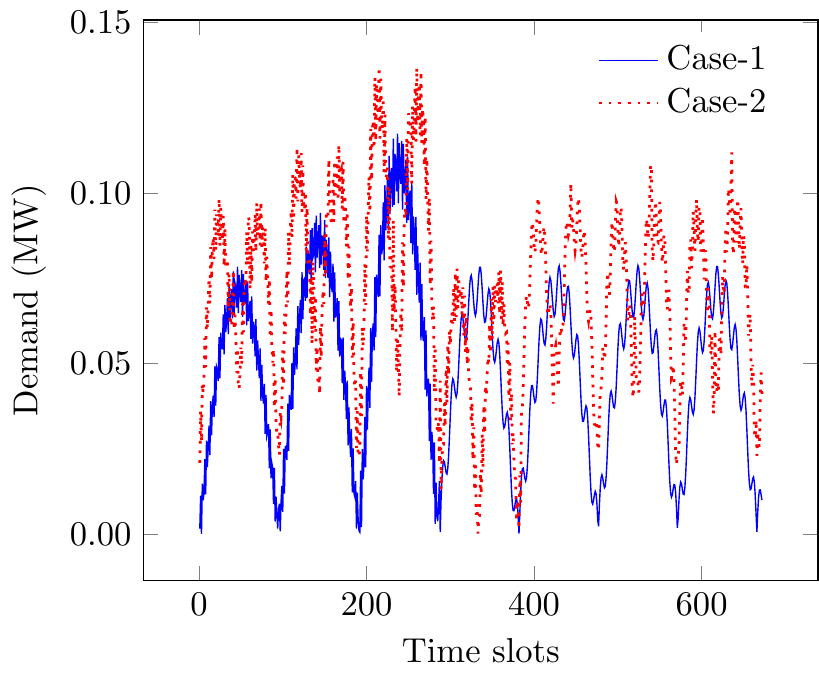}} 
   \caption{Typical load profiles on a particular bus.} 
   \label{fig:load_profile}
\end{figure}

For all the computational experiments, we use a standard IEEE 13-node radial distribution test feeder network \cite{kersting2001radial} that was modified to support a positive-sequence model. A schematic of the 13-node feeder network is shown in Fig. \ref{fig:schematic}. The algorithm is allowed to build diesel generators or batteries at the buses shown in the schematic (based on \cite{madathil2017resilient}).  
To test the effectiveness of the MPC algorithm on different loading conditions, two load profiles based on data from a New Mexico distribution utility (Kit Carson Electric Cooperative) were generated. Both high and low frequency fluctuations were utilized. The typical daily and weekly load profiles are shown in the Fig. \ref{fig:load_profile}. A time discretization of $\bm \Delta \bm t = 15$ minutes was chosen for all the runs of the algorithm. This value concurs with the time discretization of \cite{madathil2017resilient}. The load profiles are also assumed to be available at every time discretization point. 

For every battery, the charging and discharging efficiency values, $\bm \eta_{\operatorname{ch}}^b$ and $\bm \eta_{\operatorname{dis}}^b$ respectively, were set to $80\%$ and $70\%$ (from previous work \cite{felder2014optimal}). The battery's fixed cost ($\bm f_b$) and variable cost ($\bm g_b$) were \$100 and 300\$/MVA, respectively. The efficiency of diesel generators was set to $50\%$ ($\bm \eta^d$). The fixed cost and variable cost coefficients ($\bm f_d,\bm c_d^0,\bm c_d^1,\bm c_d^2$) for three different options of diesel generators were (200, 6, 35, 50), (300, 3, 10, 20), (350, 2, 5, 10).

Each run of the MPC algorithm was limited to $3$ iterations, \textit{i.e.}, $\bm N = 3$. The value of the load-shedding penalty in Eq. \eqref{objective} was set to $10^7$\$/MW. The number of stages for every run of the algorithm was set to $S = 6$. All algorithms were implemented using Julia/JuMP \cite{dunning2015jump} and Gurobi v$7.0$ was used to solve the MIQCQP. All the runs were performed on a Dell Precision T5500 workstation (Intel Xeon $E5360$ processor @$2.53$ GHz, $12$ GB RAM).

\subsection{Performance of the MPC algorithm} \label{subsec:performance}
We first examine the run times of the MPC algorithm for solving the {\problem } with different time horizons. 
The run times of the MPC algorithm are detailed in Table \ref{tab:times}.
\begin{table}[htp]
        \centering
        \caption{Run times}
        \begin{tabular}{c|c|c|c}
            \toprule
            Load & Time & Initialization & Avg. time per iteration \\
            profile & horizon & step (sec.) & MPC  (sec.)  \\
            \midrule                       
            \parbox[t]{2mm}{\multirow{3}{*}{\rotatebox[origin=c]{90}{Case 1}}} & 3 days & 67.50 & 44.98 \\ 
            & 7 days & 249.17 & 84.15 \\ 
            & 14 days &700.30 & 136.97 \\ 
            \midrule 
            \parbox[t]{2mm}{\multirow{3}{*}{\rotatebox[origin=c]{90}{Case 2}}} & 3 days & 233.79 & 31.20 \\ 
            & 7 days &313.03 &82.01 \\ 
            & 14 days & 1285.07 & 429.00 \\ 
            \bottomrule
        \end{tabular}
        \label{tab:times}
\end{table}
The average run time per iteration 
increases as the time horizon increases because the algorithm always decomposes the problem into $6$ stages irrespective of the time horizon. 
The initialization step in Table \ref{tab:times} is the time taken to compute the initial set of dual values for the first iteration of the MPC. The dual values are based on the solution to the continuous relaxation of the \problem. 

Table \ref{tab:solution} presents the solution quality measured as a relative gap from the objective value of a continuous, perspective-based relaxation of the {\problem}. The perspective-based relaxation is not presented in this article due to page limitations. For the sake of comparison with existing approaches in the literature, we also present the results obtained by using the RH algorithm. 
 \begin{table}[h!]
    \centering
    \caption{Solution quality. Load shed (LS) values are in MW.}
    \begin{tabular}{c|c|c|c|c}
        \toprule
        Load & Time & \# discretization & \multicolumn{2}{c}{Relative gap (\%)} \\
        profile & horizon & points & MPC & RH \\
        \midrule
        \parbox[t]{2mm}{\multirow{3}{*}{\rotatebox[origin=c]{90}{Case 1}}} & 3 days & 288 & 4.34 & LS (0.35)\\
        & 7 days & 672 & 3.17 & LS (2.77)\\
        & 14 days & 1344 & 3.67 & LS (0.53) \\
        \midrule
        \parbox[t]{2mm}{\multirow{3}{*}{\rotatebox[origin=c]{90}{Case 2}}} & 3 days & 288 & 4.64 & 5.04 \\
        & 7 days & 672 & 1.30 & 1.35 \\
        & 14 days & 1344 & 0.97 & 0.98\\
        \bottomrule
    \end{tabular}
    \label{tab:solution}
\end{table}
Whenever the RH algorithm sheds load, ``LS'' is reported in the relative gap column of Table \ref{tab:solution}. The actual real load shed is specified in parentheses. It is clear from the Table \ref{tab:solution} that the MPC produces solutions that are within $5$\% of the lower bound for all the cases. The RH is only able to produce comparable solutions for case-2 load profiles. These have a faster change in frequencies than case-1. 

\subsection{Battery usage} \label{subsec:battery_usage} 
\begin{figure}[!h]
   \centering
  \subfigure[Total generation cost; the value inside the bars are the real power shed.]{
  \includegraphics[scale=0.87]{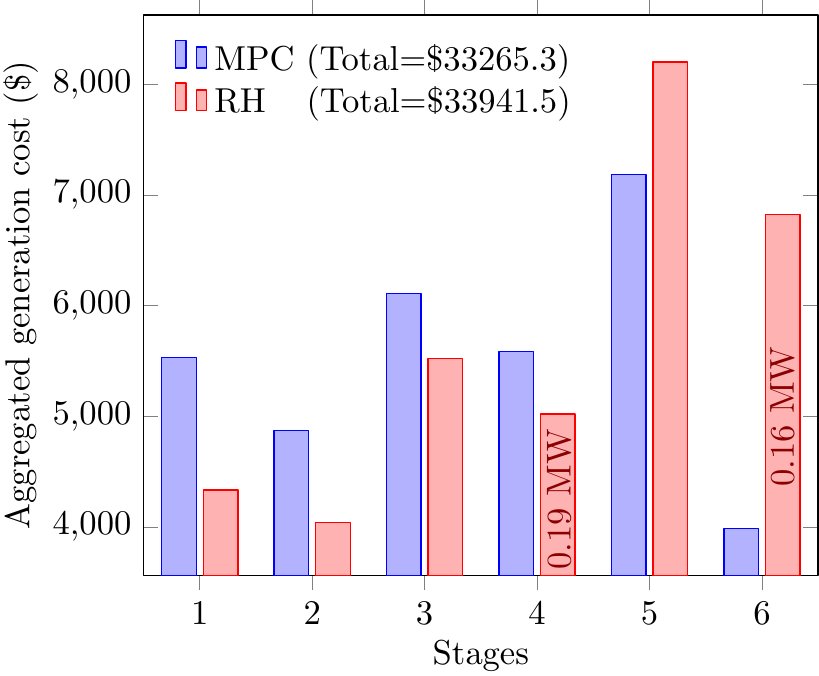}}     
   \subfigure[Aggregated state-of-charge of the batteries]{
   \includegraphics[scale=0.87]{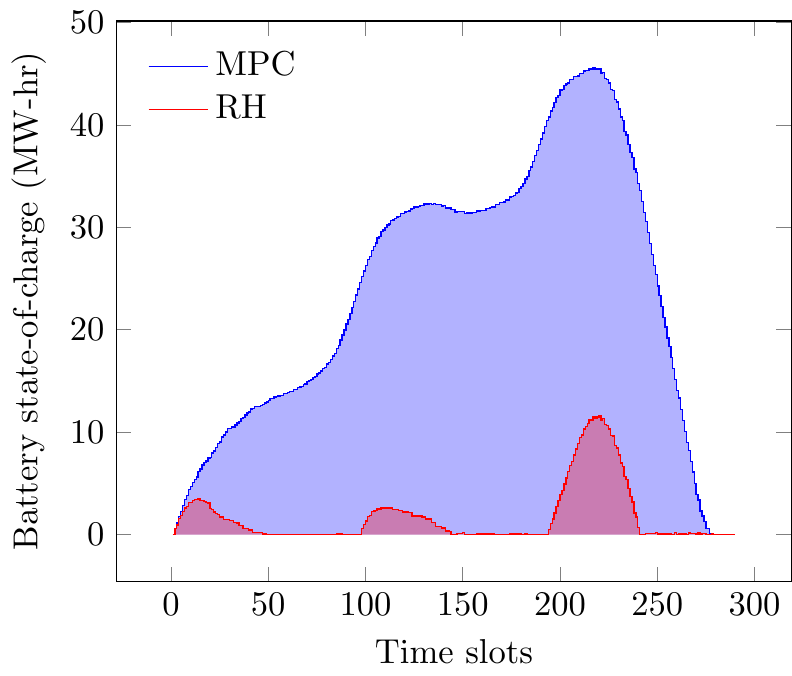}} 
   \caption{State-of-charge of batteries based on solutions with MPC vs. RH algorithms} 
   \label{fig:bat_usg}
\end{figure}

This section compares the battery usage of the MPC and RH approaches. Fig. \ref{fig:bat_usg}(b) shows the total state-of-charge of all the batteries that were built in the system by the MPC and RH solutions for the case-1 load profile (slow frequencies) with a time horizon of $3$ days. It is clearly evident from Fig. \ref{fig:bat_usg}(b) that the MPC algorithm builds and utilizes the batteries in a much more efficient way than the RH algorithm. This is attributed to the fact that the \textit{cost-to-go}-like term in the objective of the stage problems in Eq. \eqref{eq:stage_problems} better approximates the slow dynamics in the system for non-zero dual values. 
The non-zero values provide the MPC algorithm with a look-ahead feature.
The RH algorithm lacks this look-ahead feature because it is a one-pass algorithm that does not use dual values to update the solution. The short-sightedness of the RH algorithm is also evident from  Fig. \ref{fig:bat_usg}(a) which shows the total generation cost for all the stages. Despite the RH having very low generation cost in the initial few stages, the lack of a feature in the RH algorithm to foresee the load profiles in the future stages leads to a higher overall generation cost in stages $5$ and $6$ of Fig. \ref{fig:bat_usg}(a). Despite real load shedding of $0.35$MW, the aggregated generation cost of RH algorithm is 2.04\% higher than the MPC's solution. 

\subsection{Exactness of convex relaxation for battery efficiencies} \label{subsec:battery_exactness}
We remark that despite using a relaxation for modeling battery efficiencies in Sec. \ref{sec:formulation} to reduce the computational burden, the solutions obtained using the MPC algorithm for all the runs were always on the charging or discharging efficiency curves in Fig. \ref{fig:eff_modeling}. This observation is primarily due to the cost-minimization objective which ensures no further losses during charging/discharging states of the battery.  The theoretical analysis and justification is an interesting direction of future work.


\section{Conclusions}
\label{sec:conc}
In this paper, we  presented an MIQCQP formulation for long time-horizon planning and operation of a microgrid. An elegant convex relaxation to model battery efficiencies that was empirically exact was developed. We also develop a fast, scalable, hierarchical predictive control algorithm to compute feasible solutions for an MIQCQP. 
Extensive computational experiments illustrate that the algorithm scales well with increases in the time horizon and is able to compute feasible solutions that are within $5$\% of a lower bound. Future research should consider N-1 security constraints, battery degradation models, inclusion of fluctuating renewable energy sources like wind and solar, and algorithmic enhancements such as hierarchical coarse-grid-based dual updates and sliding horizon schemes with stage overlaps.

\bibliographystyle{IEEEtran}
\bibliography{references.bib}

\end{document}